\documentclass[12pt,twoside]{article}
\usepackage{amsmath,amsbsy,amsfonts}
\usepackage{graphicx,color}
\newtheorem{theorem}{Theorem}[section]
\newtheorem{lemma}[theorem]{Lemma}
\newtheorem{proposition}[theorem]{Proposition}
\newtheorem{corollary}[theorem]{Corollary}
\newtheorem{remark}[theorem]{Remark}
\newtheorem{claim}[theorem]{Claim}

\newcommand{\fim}{\hfill\rule{2mm}{2mm}}

\newcommand{\ds}{\displaystyle}
\newcommand{\R}{\mathbb{R}}

\begin{document}

\setlength{\baselineskip}{4.5mm} \setlength{\oddsidemargin}{8mm}
\setlength{\topmargin}{-3mm}
\title{\Large\sf Multiplicity of multi-bump type nodal solutions for a class of elliptic problems
with exponential critical growth in $\mathbb{R}^2$.}
\author{\sf
 Claudianor O. Alves\thanks{Partially supported by CNPq - Grant 304036/2013-7, coalves@dme.ufcg.edu.br} \;\; and \;\;  Denilson S. Pereira\thanks{denilsonsp@dme.ufcg.edu.br} \\
 Universidade Federal de Campina Grande\\
 Unidade Acad\^emica de Matem\'atica - UAMat\\
 CEP: 58.429-900 - Campina Grande - PB - Brazil\\ }

\pretolerance10000
\date{}
\numberwithin{equation}{section} \maketitle
\begin{abstract}
In this paper, we establish the existence and multiplicity of
multi-bump nodal solutions for the following class of problems
$$
-\Delta u+(\lambda V(x)+1)u=f(u),~~\mbox{in}~~\mathbb{R}^2,
$$
where $\lambda\in(0,\infty)$, $f$ is a continuous function with
exponential critical growth and $V:\mathbb{R}^2\to\R$ is a
continuous function verifying some hypotheses.

\end{abstract}


\vspace{0.3cm}


\section{Introduction}

\hspace{0.5 cm} In the present  paper, we study the existence and multiplicity of multi-bump nodal solutions for the following
class of problems
$$
\left\{
\begin{array}{l} -\Delta u+(\lambda V(x)+1)u=f(u),~~  \mbox{ in }~~ \R^2, \\[0.2cm]
u\in H^1(\R^2),
\end{array}
\right. \eqno (P)_\lambda
$$
where $\lambda\in(0,\infty)$, $V:\mathbb{R}^2\to\R$ is a nonnegative continuous
function and $f$ is a continuous function having  an exponential critical growth at
$\pm\infty$, i.e., there exists $\alpha_0>0$ such that
$$
\lim_{|t|\to +\infty}\dfrac{|f(s)|}{e^{\alpha s^2}}=0,~~\forall
\alpha>\alpha_0;~~~\lim_{|t|\to +\infty}\dfrac{|f(s)|}{e^{\alpha
s^2}}=+\infty,~~\forall \alpha<\alpha_0.
$$

There are a lot of papers concerning with existence and multiplicity
of positive solutions for $(P)_\lambda$ in $\R^N$, where the behavior
of function $V$ plays an important role. For the case $N\geq 3$, we would like to cite the papers due to Bartsch and
Wang \cite{BWa1}, Clapp and Ding \cite{MD}, Bartsch, Pankov and Wang
\cite{BPW}, Gui \cite{G}, Ding and Tanaka \cite{DT}, Alves, de
Morais Filho and Souto \cite{AMFS} and references therein.

In \cite{DT}, Ding and Tanaka considered the problem $(P)_\lambda$
assuming that $\Omega=int~ V^{-1}(\{0\})$ has $k$ connected
components and $f(s)=|s|^{q-2}s$ with $2<q<\frac{2N}{N-2}$. In that paper, it was proved that 
$(P)_\lambda$ has at least $2^k-1$ multi-bump positive
solutions for large $\lambda$. The same type of problem was
considered by Alves, de Morais Filho and Souto in \cite{AMFS} and
Alves and Souto \cite{AS}, by assuming that $f$ has a critical
growth for the case $N\geq 3$ and exponential critical growth when $N=2$, respectively.

In \cite{A1}, motivated by \cite{MD} and \cite{DT}, Alves has 
considered the existence and multiplicity of multi-bump nodal
solutions for $(P)_\lambda$, when the nonlinearity $f$ has a subcritical growth.

The motivation of our work comes from the references mentioned above, once we have observed that until moment, the existence and multiplicity of multi-bump nodal solution for $(P)_\lambda$, when $f$ has exponential critical growth in $\R^2$, were not  considered. Here, we have used a different approach in some estimates, because in our opinion, some properties that are valid for $N\geq 3$, we cannot repeat for the class of problems here studied, therefore a careful analysis is needed. 

Here, we use a result related to the existence of least energy nodal solutions for the Dirichlet Problem on a bounded domain due to Alves and Pereira \cite{AP}, which is a version of results due to  Bartsch, Weth and Willem \cite{BWW} (see also Bartsch and Weth \cite{BW}) for critical growth in $\R^2$. Moreover, we mention that the technique developed in \cite{AP} can employ to prove an existence result of least energy nodal solutions for a class of elliptic problems on a bounded domain with Neumann boundary condition. These solutions play an important role in our arguments to get multi-bump nodal solutions. Furthermore, as in \cite{A1}, we modify all the sets that appear in the minimax arguments found in \cite{DT} to get nodal solutions. Our main result completes the studies made in \cite{DT}, \cite{AS} and \cite{A1}, in the following points:
\begin{itemize}
\item In \cite{DT}, the nonlinearity is homogeneous with subcritical growth, the solutions are positive  and $N \geq 3$.

\item In \cite{AS}, the nonlinearity has an exponential critical growth and $N=2$, but the solutions  are positive.

\item In \cite{A1}, the solution are multi-bump nodal, however the nonlinearity has subcritical growth in $\R^N$ and $N\geq 3$.
\end{itemize}

We would like to mention that problems involving exponential critical growth have received a special attention at last years, see for example, \cite{AdoOM,Cao,doORuf,DMR,DdOR} for semilinear elliptic equations, and \cite{1,2,3,4,5} for quasilinear equations.,

Since we will work with exponential critical growth, some versions of the Trudinger-Moser inequality are crucial in arguments. The first version that we would like to recall is due to Trundiger and Moser, see \cite{M} and \cite{T}, which claims if $\Omega$ is a bounded domain with smooth boundary, then for any  $ u \in H^{1}_0(\Omega)$,
\begin{equation} \label{X0}
\int_{\Omega}
e^{\alpha u^2}dx
< +\infty, \,\,\,\, \mbox{ for every }\,\,\alpha >0.
\end{equation}
Moreover, there exists a positive constant $C=C(\alpha,|\Omega|)$ such that
\begin{equation} \label{X1}
\sup_{||u||_{H_0^{1}(\Omega)} \leq 1} \int_{\Omega} e^{\alpha u^2} dx \leq C , \,\,\,\,\,\,\, \forall \, \alpha  \leq 4 \pi .
\end{equation}

A version in $H^{1}(\Omega)$ has been proved by Adimurthi and Yadava \cite{AY}, and it says that if $\Omega$ is again a bounded domain with smooth boundary, then for any  $ u \in H^{1}(\Omega)$,
\begin{equation} \label{X02}
\int_{\Omega}
e^{\alpha u^2}dx
< +\infty, \,\,\,\, \mbox{ for every }\,\,\alpha >0.
\end{equation}
Furthermore, there exists a positive constant $C=C(\alpha,|\Omega|)$ such that
\begin{equation} \label{X11}
\sup_{||u||_{H^{1}(\Omega)} \leq 1} \int_{\Omega} e^{\alpha u^2} dx \leq C , \,\,\,\,\,\,\, \forall \, \alpha  \leq 2 \pi .
\end{equation}

The third version that we will use is due to Cao \cite{Cao}, which is version of the Trundiger-Moser inequality in whole space $\R^2$ and has the following statement:
\begin{equation} \label{X2}
\int_{\R^2}\left(e^{\alpha u^2}-1\right)dx<+\infty,\ \ \mbox{for all} \,\,\, u\in H^1(\R^2)\ \ \mbox{and} \,\, \alpha>0.
\end{equation}
Besides, if $\alpha < 4\pi$ and $|u|_{L^2(\R^2)}\leq M$, there exists
a constant \linebreak $C_1=C_1(M,\alpha)$ verifying
\begin{equation} \label{X4}
\sup_{ |\nabla u|_{L^2(\R^2)}\leq 1}
\int_{\R^2}\left(e^{\alpha u^2}-1\right)dx\leq C_1.
\end{equation}

In what follows, let us denote $\Omega=int~V^{-1}(\{0 \})$ and we
suppose that
\begin{enumerate}
\item[$(H1)$]$\Omega$ is non-empty, bounded, $\partial\Omega$
is smooth and $V^{-1}(\{0 \})=\overline{\Omega}$;

\item[$(H2)$] $\Omega$ has $k$ connected components denoted by
$\Omega_j$, $j\in\{1,...,k\}$, such that $dist(\Omega_i,\Omega_j)>0$ for $i\neq j$.
\end{enumerate}

Hereafter, the function $f$ satisfies the ensuing assumptions:

\begin{enumerate}

\item[$(f_1)$] There is $C>0$ such that
$$
|f(s)|\leq Ce^{4\pi |s|^2}\ \ \mbox{for all}\ \ s\in\R;
$$

\item[$(f_2)$] $\ds\lim_{s\rightarrow
0}\dfrac{f(s)}{s}=0$;

\item[$(f_3)$] There is $\theta>2$ such that
$$0<\theta F(s):=\theta\int_0^{s}f(t)dt\leq sf(s),\ \ \mbox{for all}\ \
s\in\R\setminus\{0\}.$$

\item[$(f_4)$] The function $s\rightarrow\dfrac{f(s)}{|s|}$ is
strictly increasing in $(0,+\infty)$.

\item[$(f_5)$] There exist constants $p>2$ and $C_p>0$ such  that
$$
sign(s)f(s)\geq C_p|s|^{p-1}\ \ \mbox{for all}\ \
s\in\R
$$
with
$$
C_p>\left[\dfrac{4k\theta}{\theta-2}S_p\right]^{(p-2)/2},
$$
where $S_p=\displaystyle\max_{1\leq j\leq k}\gamma_j,~~\gamma_j=\displaystyle\inf_{u\in\mathcal{M}_{\Omega_j}}\phi_j(u)$,
$$
\mathcal{M}_{\Omega_j}=\left\{u\in H_0^1(\Omega_j):~ u^{\pm}\neq 0~
\mbox{and}~ \phi'_j(u^{\pm})u^{\pm}=0\right\}
$$
and
$$
\phi_j(u)=\dfrac{1}{2}\!\int_{\Omega_j}\!\!\!\left(|\nabla
u|^2+|u|^2\right)-\dfrac{1}{p}\!\int_{\Omega_j}\!\!\!|u|^p.
$$

\end{enumerate}

It is easily seen that $(f_1)-(f_5)$ hold for nonlinearities of the form
$$
f(s)=2\alpha s \left(e^{\alpha s^2}-1\right),~\mbox{for}~\alpha\in(0,4\pi).
$$

Our main  result is the following.

\noindent  \begin{theorem}\label{bd} Assume that $(f_1)-(f_5)$ and
$(H_1)-(H_2)$  hold. Then, for any non-empty subset $\Gamma$ of
$\{1,...,k\}$, there exists $\lambda^{*}>0$ such that, for
$\lambda\geq\lambda^{*}$, problem $(P)_\lambda$ has a nodal solution
$u_\lambda$. Moreover, the family
$\{u_\lambda\}_{\lambda\geq\lambda^{*}}$ has the following property:
For any sequence $\lambda_n\rightarrow \infty$, we can extract a
subsequence $\lambda_{n_i}$ such that $u_{\lambda_{n_i}}$ converges
strongly in $H^1(\R^2)$ to a function $u$ which satisfies $u(x)=0$
for $x\notin \Omega_{\Gamma}=\cup_{j \in \Gamma}\Omega_j$, and the restriction $u|_{\Omega_j}$
is a nodal solution with least energy of
$$
 -\Delta u+u=f(u),~~  \mbox{ in }~~ \Omega_j, ~~
 u|_{\partial\Omega_j}=0~~\mbox{for}~~j\in\Gamma.
$$
\end{theorem}

\vspace{0.5 cm}

The plan of the paper is as follows: In Section 2, we prove some technical results involving bounded domains, which will be useful in the proof the Theorem \ref{bd}. In Sections 3 and 4, we consider an auxiliary problem and study some properties of the energy functional associated with that problem. Finally in Section 5, we prove the main result.

\section{Preliminaries}
Throughout this paper, we will use the following notations:

\begin{itemize}
\item If $h$ is a measurable function, we denote by $\int_{\R^2}h$
the following integral $\int_{\R^2}h(x)dx$.

\item The symbols $\|u\|,~|u|_r~(r>1)$ and $|u|_\infty$ denote the
usual norms in the spaces $H^1(\R^2),~L^r(\R^2)$ and
$L^\infty(\R^2)$, respectively.

\item For an open set $\Theta\subset\R^2$, the symbols
$\|u\|_\Theta,~|u|_{r,\Theta}~(r>1)$ and $|u|_{\infty,\Theta}$
denote the usual norms in the spaces $H^1(\Theta),~L^r(\Theta)$ and
$L^\infty(\Theta)$, respectively.
\end{itemize}
From now on, we will work with the space $\mathcal{H}_\lambda$ defined
by
$$
\mathcal{H}_\lambda=\left\{u\in
H^1(\R^2):~\int_{\R^2}V(x)u^2<\infty\right\}
$$
endowed with the norm
$$
\|u\|_\lambda=\left[\int_{\R^2}\left(|\nabla u|^2+(\lambda
V(x)+1)u^2\right)\right]^{1/2}.
$$
It is easy to see that $(\mathcal{H}_\lambda,\|\cdot\|_\lambda)$ is
a Hilbert space for $\lambda>0$.

We also write for an open set $\Theta\subset\R^2$
$$
\mathcal{H}_\lambda(\Theta)=\left\{u\in H^1(\Theta):~\int_{\Theta}
V(x)u^2<\infty\right\}
$$
and
$$
\|u\|_{\lambda,\Theta}=\left[\int_{\Theta}\left(|\nabla
u|^2+(\lambda V(x)+1)u^2\right)\right]^{1/2}.
$$

As a consequence of the above considerations, there exist $\nu_0,
\delta_0>0$ with $1\approx \delta_0<1$ and $\nu_0\approx 0$ such
that for all open set $\Theta\subset\R^2$
\begin{equation}\label{1}
\delta_0\|u\|^2_{\lambda,\Theta}\leq\|u\|^2_{\lambda,\Theta}-\nu_0|u|^2_{2,\Theta},~~\forall
u\in \mathcal{H}_\lambda(\Theta)~\mbox{and}~\lambda> 0.
\end{equation}

From assumptions $(f_1)$ and $(f_2)$, given $\epsilon\!>\!0$, $q\!\geq\! 1$
and $\tau\!>\!1$, there exists a constant $C=C(\epsilon,q,\alpha)>0$
such that
\begin{equation}\label{eg}
|sf(s)|,\ \ |F(s)| \leq \epsilon s^2+C|s|^{q}b_\tau(s),\ \ \mbox{for
all}\ \ s\in \R,
\end{equation}
where
$$
b_\tau(s):=\left(e^{4\pi\tau s^2}-1\right).
$$

The below result is a consequence of Trundinger-Moser inequality given in (\ref{X4}) and its proof can be found in \cite{AS}.
\begin{corollary}\label{ctec}
Let $(u_\lambda)$ be a family in $H^1(\R^2)$ satisfying $\displaystyle\sup_{\lambda\geq 1}\|u_\lambda\|^2\leq m<1.$
For $\tau,q>1$ satisfying $\tau qm<1$, there exists
$C=C(\tau,q,m)>0$ such that $b_\tau(u_\lambda)=\left(e^{4\pi \tau
u_\lambda^2}-1\right)$ belongs to $L^q(\R^2)$ and
$$
\sup_{\lambda\geq 1}\{|b_\tau(u_\lambda)|_q\}<\infty.
$$
\end{corollary}

\subsection{Neumann and Dirichlet problems}\label{snd}

\hspace{0.5 cm} In this section, we denote by $I_j:H_0^1(\Omega_j)\to\R$ and
$\Phi_{\lambda,j}: H^1(\Omega_j')\to\R$ the following energy
functionals
$$
I_j(u)=\dfrac{1}{2}\int_{\Omega_j}(|\nabla
u|^2+u^2)-\int_{\Omega_j}F(u)
$$
and
$$
\Phi_{\lambda,j}(u)=\dfrac{1}{2}\int_{\Omega_j'}(|\nabla
u|^2+(\lambda V(x)+1)u^2)-\int_{\Omega_j'}F(u).
$$
It is well known that $I_j$ and $\Phi_{\lambda,j}$ are $C^1$ and their
critical points are weak solutions of the problems

\begin{equation}\label{2}
\left\{
\begin{array}{rcl} -\Delta u+u&\!\!=\!\!&f(u),~  \mbox{ in }~ \Omega_j, \\[0.2cm]
u&\!\!=\!\!&0,~ \mbox{on}~\partial\Omega_j
\end{array}
\right.
\end{equation}
and
\begin{equation}\label{3}
\left\{
\begin{array}{rcl} -\Delta u+(\lambda V(x)+1)u&\!\!=\!\!&f(u),~\mbox{ in }~\Omega_j', \\[0.2cm]
\dfrac{\partial u}{\partial\nu}&\!\!=\!\!& 0,~ \mbox{on}~\partial\Omega_j',
\end{array}
\right.
\end{equation}
respectively. 

Hereafter, $d_j$ and $d_{\lambda,j}$ denote the real numbers given by 
$$
d_j=\inf_{\mathcal{M}_j}I_j~~~\mbox{and}~~~d_{\lambda,j}=\inf_{\mathcal{M}_{\lambda,j}}\Phi_{\lambda,j},
$$
where $\mathcal{M}_j$ and $\mathcal{M}_{\lambda,j}$ denote the nodal Nehari sets 
$$
\mathcal{M}_j=\{u\in H_0^1(\Omega_j):~u^{\pm}\neq 0~~\mbox{and}~~I_j'(u^{\pm})u^{\pm}=0\},
$$
and
$$
\mathcal{M}_{\lambda,j}=\{u\in H^1(\Omega_j'):~u^{\pm}\neq 0~~\mbox{and}~~\Phi_{\lambda,j}'(u^{\pm})u^{\pm}=0\}.
$$

By a result found in \cite{AP}, we know that there is $w_j\in \mathcal{M}_j$ verifying
$$
I_j(w_j)=d_j~~~\mbox{and}~~~I'_j(w_j)=0.
$$
Here, we would like to point that the same approach can be employed to show that there is $w_{\lambda,j}\in \mathcal{M}_{\lambda,j}$ satisfying
$$
\Phi_{\lambda,j}(w_{\lambda,j})=d_{\lambda,j}~~~\mbox{and}~~~\Phi'_{\lambda,j}(w_{\lambda,j})=0.
$$
To see why, it remains to observe that $(f_5)$ yields if $(v_n)$ is a Palais-Smale sequence associated to $\Phi_{\lambda,j}$ at $d_{\lambda,j}$, then
$$
\limsup_{n \to +\infty}\|v_n\|_{\Omega_j'}^{2} < 1/2.
$$
The above estimate is the key point to apply the Trudinger-Moser inequality due to Adimurthi and Yadava, see  (\ref{X11}). In doing so, the reader will see that the existence of $w_{\lambda,j}$ follows as in \cite{AP}, replacing $H^{1}_0(\Omega)$ by $H^{1}(\Omega)$.


\section{An auxiliary problem}

\hspace{0.5 cm} In this section, as in \cite{A1}, \cite{PF} and \cite{DT}, we will modify conveniently the function $f$. 

To this end, let $\nu_0$ be the constant given in $(\ref{1})$, $a>0$ verifying
$\max\{f(a)/a,f(-a)/(-a)\}<\nu_0$ and $\tilde{f},\tilde{F}:\R\rightarrow \R$ the following functions
$$
\tilde{f}(s)=\left\{
\begin{array}{ll}
\dfrac{-f(-a)}{a}s & \mbox{if}~ s<-a,\\[0.4cm]
f(s)               &\mbox{if}~ |s|\leq a,\\[0.4cm]
\dfrac{f(a)}{a}s   & \mbox{if}~ s>a
\end{array}
\right.
$$
and
$$
\tilde{F}(s)=\displaystyle\int_{0}^{s}\tilde{f}(\tau)d\tau,
$$
which fulfills the inequalities
\begin{equation}\label{vc11}
\tilde{f}(s)\leq\nu_0|s|,~~~\forall s\in\R,
\end{equation}
\begin{equation}\label{vc12}
\tilde{f}(s)s\leq\nu_0|s|^2,~~~\forall s\in\R,
\end{equation}
and
\begin{equation}\label{vc12}
\tilde{F}(s)\leq\dfrac{\nu_0}{2}|s|^2,~~~\forall s\in\R.
\end{equation}

From now on, for each subset $\Gamma\subset\{1,...,k\}$, let us
consider
$$
\chi_{\Gamma}=\left\{
\begin{array}{l}
1,~~\mbox{for}~~x\in\Omega_{\Gamma}',\\[0.1cm]
0,~~\mbox{for}~~x\in\R^2\setminus \Omega_{\Gamma}',
\end{array}
\right.
$$
where $\Omega_{\Gamma}'=\displaystyle \cup_{j \in \Gamma}\Omega_j'$. Using the above functions, we define 
$$
g(x,s)=\chi_\Gamma(x)f(s)+(1-\chi_\Gamma(x))\tilde{f}(s)
$$
and
$$
G(x,s)=\int_0^sg(x,t)dt=\chi_\Gamma(x)F(s)+(1-\chi_\Gamma(x))\tilde{F}(s).
$$
It is easy to see that $g$ satisfies $(\ref{eg})$ uniformly in $x\in \R^2$, that is,
\begin{equation}\label{g}
|g(x,s)|\leq \epsilon|s|+C|s|^{q-1}b_\tau(s),~~~\forall
s\in\R,~~x\in\R^2.
\end{equation}
Using the above estimate, it follows that $\Phi_\lambda:\mathcal{H}_\lambda\rightarrow\R$ given by
$$
\Phi_\lambda(u)=\dfrac{1}{2}\int_{\R^2}(|\nabla u|^2+(\lambda
V(x)+1)u^2)-\int_{\R^2}G(x,u)
$$
belongs to $C^1(\mathcal{H}_\lambda,\R)$ and its critical points are
weak solutions of
$$
-\Delta u+(\lambda
V(x)+1)u=g(x,u)~~\mbox{in}~~\R^2.~~~~~~~~~~(A)_\lambda
$$

\begin{remark}\label{rea}
In this moment, we would like to detach that some nodal solutions of
$(A)_\lambda$ are solutions of the original problem $(P)_\lambda$.
More precisely, if $u_\lambda$ is a nodal solution of $(A)_\lambda$
verifying $|u(x)|\leq a$ in $\R^2\setminus\Omega_\Gamma'$, then it
is a nodal solution for $(P)_\lambda$.
\end{remark}

In the sequel, we study the convergence of Palais-Smale sequences
related to the functional $\Phi_\lambda$. The first of them is related to  boundedness of these sequences. However, it follows repeating the same arguments explored in \cite[Lemma 3.1]{AS}, then we will omit its proof. 


\begin{lemma}\label{l1}
If $(u_n)$ is a $(PS)_c$ sequence to $\Phi_\lambda$, then
$$
\limsup_{n\rightarrow \infty}\|u_n\|_\lambda^2\leq\dfrac{2\theta
c}{\delta_0(\theta-2)},
$$
where $\delta_0$ is given in $(\ref{1})$.
\end{lemma}

\vspace{0.5 cm}

In the next, we denote by $D$ the ensuing real number
$$
D=\sum_{j=1}^kd_j.
$$
This number is very special for us, because we will show that $\Phi_\lambda$ verifies the well known Palais-Smale in $(0,D]$. To prove this fact, we need of the following estimate from above for $D$.

\begin{lemma}\label{l2}
If $(f_1)-(f_5)$ holds, then
$0<D<\dfrac{\delta_0(\theta-2)}{4\theta}$.
\end{lemma}

\noindent {\bf Proof.} In order to prove this inequality, for each
$j\in\{1,...,k\}$, we fix a nodal function $v_j\in H_0^1(\Omega_j)$
such that $v_j\in \mathcal{M}_{\Omega_j}$ and
\begin{equation}\label{5}
\phi_j(v_j)=\gamma_j.
\end{equation}
The reader can find the proof of the existence of such functions in \cite{BW}. Since $v_j^{\pm}\neq 0$, there exist
$s_j,~t_j>0$ such that $s_jv_j^{+}+t_jv_j^{-}\in\mathcal{M}_j$.
Thus,
$$
d_j\leq I_j(s_jv_j^{+}+t_jv_j^{-})=I_j(s_jv_j^{+})+I_j(t_jv_j^{-}),
$$
or equivalently,
$$
\begin{array}{ll}
d_j\leq&\dfrac{s_j^2}{2}\displaystyle\int_{\Omega_j}(|\nabla
v_j^{+}|^2+|v_j^{+}|^2)-\displaystyle\int_{\Omega_j}F(s_jv_j^{+})\\[0.3cm]
&+\dfrac{t_j^2}{2}\displaystyle\int_{\Omega_j}(|\nabla
v_j^{-}|^2+|v_j^{-}|^2)-\displaystyle\int_{\Omega_j}F(t_jv_j^{-}).
\end{array}
$$
Using the fact that $v_j^{\pm}\in \mathcal{M}_{\Omega_j}$ and $(f_4)$, we obtain
$$
d_j\leq\left\{\dfrac{s_j^2}{2}-\dfrac{C_ps_j^p}{p}\right\}\int_{\Omega_j}|v_j^{+}|^p+\left\{\dfrac{t_j^2}{2}-\dfrac{C_pt_j^p}{p}\right\}\int_{\Omega_j}|v_j^{-}|^p.
$$
Then,
$$
d_j\leq\max_{r\geq
0}\left\{\dfrac{r^2}{2}-\dfrac{C_pr^p}{p}\right\}\int_{\Omega_j}|v_j|^p.
$$
Noting that
$$
\max_{r\geq
0}\left\{\dfrac{r^2}{2}-\dfrac{C_pr^p}{p}\right\}=C_p^{\frac{2}{2-p}}\left(\dfrac{1}{2}-\dfrac{1}{p}\right),
$$
it follows
$$
d_j\leq
C_p^{\frac{2}{2-p}}\left(\dfrac{1}{2}-\dfrac{1}{p}\right)\int_{\Omega}|v_j|^p.
$$
Combining $(\ref{5})$ with the above inequality, we derive
$$
d_j\leq C_p^{\frac{2}{2-p}}\gamma_j,
$$
and so,
$$
D=\sum_{j=1}^kd_j\leq k S_p C_p^{\frac{2}{2-p}}<
\dfrac{\theta-2}{4\theta}.
$$
Since $\delta_0$ can be chosen close to $1$, the last inequality
leads to 
$$
D<\dfrac{\delta_0(\theta-2)}{4\theta}.
$$
\fim

\begin{proposition}\label{psc}
For $\lambda\geq 1$, the functional $\Phi_\lambda$ satisfies
$(PS)_c$ condition for all $c\in (0,D]$. More precisely, any
$(PS)_c$ sequence $(u_n)\subset\mathcal{H}_\lambda$ to
$\Phi_\lambda$ has a strongly convergent subsequence in
$\mathcal{H}_\lambda$.
\end{proposition}

\noindent {\bf Proof.} \, Let $(u_n)\subset\mathcal{H}_\lambda$ be a Palais-Smale sequence for $\Phi_\lambda$ at the level $c\in(0,D]$. By
Lemmas \ref{l1} and \ref{l2},
$$
\limsup_{n\to\infty}\|u_n\|_{\lambda}^2\leq \dfrac{2\theta D}{\delta_0(\theta-2)}< \dfrac{1}{2}.
$$
Thus, $(u_n)$ is a bounded sequence in $\mathcal{H}_\lambda$. Since $\mathcal{H}_\lambda$ is a reflexive Banach space, there exists $u\in\mathcal{H}_\lambda$ such that, for some subsequence, still denoted by $(u_n)$,
$$
u_n\rightharpoonup u~~\mbox{in}~~\mathcal{H_\lambda}, \,\,\,\,\, u_n\rightharpoonup u~~\mbox{in}~~H^1(\R^2)~~~~\mbox{and}~~~~ u_n\rightarrow u~~\mbox{in}~~ L^s_{loc}(\R^2),~~\forall s\geq 1.
$$
Using similar arguments as in \cite[Lemma 1.1]{PF}, for each $\epsilon>0$, there exists $R>0$ such that
\begin{equation}\label{6}
\limsup_{n\to\infty}\int_{\R^2\setminus B_R(0)}\left(|\nabla u_n|^2+(\lambda V(x)+1)|u_n|^2\right)\leq
\epsilon.
\end{equation}

\begin{claim}\label{beg1}
The following limits occur
\begin{enumerate}
\item[(a)] $\displaystyle \int_{\R^2}g(x,u_n)u_n\to\int_{\R^2}g(x,u)u$;
\item[(b)] $\displaystyle \int_{\R^2}g(x,u_n)v\to\int_{\R^2}g(x,u)v,~~~~\forall v\in \mathcal{H}_\lambda$.
\end{enumerate}
\end{claim}
In fact, from $(\ref{g})$,
$$
|g(x,u_n)u_n|\leq \eta|u_n|^2+C_\eta|u_n|b_\tau(u_n),~~\forall x\in \R^2,~n\in\mathbb{N}.
$$
Then, setting the functions
$$
h_n:=\eta|u_n|^2+C_\eta|u_n|b_\tau(u_n)~~~\mbox{and}~~~ h:=\eta|u|^2+C_\eta|u|b_\tau(u),
$$
it follows that $|g(x,u_n)u_n|\leq h_n(x)$. Since, $u_n\to u$ in $L^s_{loc}(\R^2)$, $\forall s\geq 1$, we can assume that 
$$
u_n(x)\to u(x)~~~~\mbox{a.e. in}~~B_R(0), 
$$
for any $R>0$. Then,
$$
g(x,u_n(x))u_n(x)\to g(x,u(x))u(x)~~~~\mbox{a.e. in}~~B_R(0)
$$
and
$$
h_n(x)\to h(x)~~~~\mbox{a.e. in}~~B_R(0).
$$
We claim that
$$
h_n\to h~~~~\mbox{em}~~L^1(B_R(0)).
$$
Indeed, since $\displaystyle\limsup_{n\to\infty}\|u_n\|_\lambda^2< 1$, for a fixed $m\in(0,1)$, there is a subsequence of $(u_n)$, still denoted by $(u_n)$, such that
$$
\sup_{n\geq 1}\|u_n\|^2\leq m<1.
$$
Fixing $q, \tau>1$ sufficiently close to $1$, such that $\tau qm<1$, by Corollary \ref{ctec} there exists
$C>0$ such that $b_\tau(u_n)\in L^q(\R^2)$ with
$$
|b_\tau(u_n)|_q\leq C,~~\forall n\in\mathbb{N}.
$$
Thereby, $(b_\tau(u_n))_n$ is a bounded sequence in $L^q(B_R(0))$ and  
$$
b_\tau(u_n)\rightharpoonup b_\tau(u)~~~\mbox{in}~~ L^q(B_R(0)).
$$
Since
$$
|u_n|\to|u|~~~\mbox{em}~~ L^{q'}(B_R(0)),~~\mbox{where}~~1/q+1/q'=1,
$$
we have
$$
|u_n|b_\tau(u_n)\to |u|b_\tau(u)~~~\mbox{in}~~ L^1(B_R(0)).
$$
From this,  $h_n\to h$ in $L^1(B_R(0))$, and so, 
$$
\lim_{n\to\infty}\left|\int_{B_R(0)}g(x,u_n)u_n-\int_{B_R(0)}g(x,u)u\right|=0,
$$
for any $R>0$. Now, in what follows, we set 
$$
L_{n,1}:=\int_{\R^2\setminus B_R(0)}|g(x,u_n)u_n-g(x,u)u|.
$$
From (\ref{vc12}),
$$
|g(x,t)t|=\tilde{f}(t)t\leq \nu_0|t|^2,~~\forall x\in\R^2\setminus B_R(0),~t\in\R.
$$
Thus, 
$$
\begin{array}{ll}
L_{n,1}&\leq\displaystyle\int_{\R^2\setminus B_R(0)}\nu_0|u_n|^2+\nu_0\displaystyle\int_{\R^2\setminus B_R(0)}|u|^2\\[0.3cm]
&\leq \displaystyle\int_{\R^2\setminus B_R(0)}\left(|\nabla u_n|^2+(\lambda V(x)+1)|u_n|^2\right)+\nu_0\displaystyle\int_{\R^2\setminus B_R(0)}|u|^2.
\end{array}
$$
From (\ref{6}), given $\epsilon >0$, there is $R>0$ such that 
$$
\limsup_{n \to +\infty}\int_{\R^2\setminus B_R(0)}\left(|\nabla u_n|^2+(\lambda V(x)+1)|u_n|^2\right) \leq \epsilon
$$
Since $u\in L^2(\R^2)$, increasing $R$ if necessary, we also can suppose that 
$$
\int_{\R^2\setminus B_R(0)}|u|^2\leq\dfrac{\epsilon}{\nu_0}.
$$
In doing so, we get
$$
\limsup_{n\to\infty}L_{n,1}\leq 2\epsilon,~~\forall \epsilon>0,
$$
implying that
$$
\lim_{n\to\infty} L_{n,1}=0.
$$
From this, we have $(a)$. The proof of $(b)$ follows using the same argument.

Now, recalling that
$$
\|u_n-u\|_{\lambda}^2=\|u_n\|_{\lambda_n}^2-2\langle u_n,u\rangle_{\lambda}+\|u\|_{\lambda}^2,
$$
the limits $\Phi'_\lambda(u_n)u_n=o_n(1)$ and $\Phi'_\lambda(u_n)u=o_n(1)$ lead to 
\begin{equation}\label{moz1}
\|u_n-u\|_\lambda^2=\int_{\R^2}g(x,u_n)u_n-\int_{\R^2}g(x,u_n)u+o_n(1).
\end{equation}
Combining the last equality with the Claim \ref{beg1}, we derive
$$
u_n\to u~~~\mbox{in}~~\mathcal{H}_\lambda~~ \mbox{and}~~ H^1(\R^2),
$$
showing that $\Phi_\lambda$ satisfies the $(PS)_c$ condition, for $c\in(0,D]$.  \hspace{1 cm} \mbox{} \fim

\vspace{0.5 cm}

Our next goal is to study the behavior of a generalized Palais-Smale
sequence corresponding to a sequence of functionals. From now on, we
say that $(u_n)\subset H^1(\R^2)$ is a $(PS)_{\infty,c}$ sequence,
if there exist $\lambda_n\rightarrow\infty$ such that
$u_n\in\mathcal{H}_{\lambda_n}$ verifying 
$$
\Phi_{\lambda_n}(u_n)\rightarrow c~~\mbox{and}~~
\|\Phi_{\lambda_n}'(u_n)\|_{\lambda_n}^{*}\rightarrow
0.~~~~~~~~~~(PS)_{\infty,c}
$$
The proof of the next proposition follows with the same arguments found in  \cite[Proposition 3.2]{AS}, then we will omits it proof. 

\begin{proposition}\label{p1}
Let $(u_n)$ be a $(PS)_{\infty,c}$ sequence with $c\in  (0,D]$.
Then, for some subsequence, still denoted by $(u_n)$, there exists
$u\in H^1(\R^2)$ such that
$$u_n\rightarrow u~~\mbox{in}~~ H^1(\R^2).$$
Moreover,
\begin{enumerate}

\item[$(i)$] $u\equiv 0$ in $\R^2\setminus\Omega_\Gamma$ and
$u|_\Omega$ is a solution of
$$
\left\{
\begin{array}{l} -\Delta u+u=f(u),~~  \mbox{ in }~~ \Omega_j, \\
\mbox{}\\
u=0,~~on~~\partial\Omega_j,
\end{array}
\right. \eqno (P)_j
$$
for each $j\in\Gamma$;
\item[$(ii)$] $\|u_n-u\|_{\lambda_n}\rightarrow 0$;

\item[$(iii)$] 
$$
\lambda_n\int_{\R^2}V(x)|u_n|^2\rightarrow
0,~~~\|u_n\|_{\lambda_n,\R^2\setminus\Omega_\Gamma}^2\rightarrow 0
$$
and
$$
\|u_n\|_{\R^2\setminus\Omega_j'}^2\rightarrow\int_{\Omega_j}(|\nabla
u|^2+u^2)~~~\mbox{for all}~~~j\in\Gamma.
$$
\end{enumerate}
\end{proposition}

Now, we are able to study the boundedness outside $\Omega_\Gamma'$ for some solutions of $(A)_\lambda$. To this end, we will use the Moser iteration technique \cite{M1}, adapting arguments found in \cite{AG} and \cite{AS}.
\begin{proposition}\label{p500}
Let $\{ u_\lambda\}\subset\mathcal{H}_\lambda$ be a family of nodal solution of $(A)_\lambda$
with $\|u_\lambda\|^2\leq m<1$ for all $\lambda\geq 1$. Then, there exists $K>0$ such that
$$
|u_\lambda|_{\infty}\leq K,~~~\forall\lambda\geq 1.
$$
\end{proposition}
\noindent{\bf Proof.} The basic idea is as following: For each $\lambda\!\geq\! 1$, $L\!>\!0$ and $\beta\!>\!1$, let
$$
u_{L,\lambda}^{+}:=\left\{
\begin{array}{lll}
u_{\lambda}^{+},& \mbox{if}~~~u_\lambda\leq L,\\[0.3cm]
L,&\mbox{if}~~~u_\lambda\geq L,
\end{array}
\right.
$$
$$
z_{L,\lambda}^{+}:=(u_{L,\lambda}^{+})^{2(\beta-1)}u_\lambda^{+}~~~~~~ \mbox{and}~~~~~~ w_{L,\lambda}^{+}:=u_\lambda^{+}(u_{L,\lambda}^{+})^{\beta-1}.
$$
Using the fact that $u_\lambda$ is a nodal solution to $(A_\lambda)$ and taking $z_{L,n}^{+}$ as a test function, we obtain
\begin{equation}\label{u1}
\int_{\R^2}\nabla u_\lambda^{+}\nabla z_{L,\lambda}^{+}+\int_{\R^2}(\lambda V(x)+1)u_\lambda^{+} z_{L,\lambda}^{+}=\int_{\R^2}g(x,u_\lambda^{+})z_{L,\lambda}^{+}.
\end{equation}
Recalling that given $\epsilon >0$ there is $C_\epsilon>0$ such that
\begin{equation}\label{u2}
g(x,u_\lambda^{+})\leq \epsilon u_\lambda^{+}+C_\epsilon b_\tau(u_\lambda^{+})u_\lambda^{+},
\end{equation}
where $b_\tau(u_\lambda^{+})\in L^q(\R^2)$ for some $q>1$, $q\approx 1$ with
\begin{equation}\label{u3}
|b_\tau(u_\lambda^{+})|_q\leq C,~~~~ \forall \lambda\geq 1,
\end{equation}
it follows from $(\ref{u2})$ and $(\ref{u1})$, 
$$
|w_{L,\lambda}^{+}|_\gamma^2 \leq C\beta^2\int_{\R^2}b_\tau(u_\lambda^{+})|w_{L,\lambda}^{+}|^2.
$$
Using H\"older's inequality,
$$
|w_{L,\lambda}^{+}|_\gamma^2 \leq C\beta^2\left[\int_{\R^2}|b_\tau(u_\lambda^{+})|^q\right]^{1/q}\left[\int_{\R^2}|w_{L,\lambda}^{+}|^{2q'}\right]^{1/q'},
$$
where $1/q+1/q'=1$. Then, by $(\ref{u3})$,
$$
|w_{L,\lambda}^{+}|_\gamma^2 \leq C\beta^2|w_{L,\lambda}^{+}|_{2q'}^{2},~~~~\forall \lambda\geq 1,
$$
for any $L>0$, $\beta>1$ and $\gamma\geq 2$, where $C>0$ depends only on $\gamma$.

Note that by Sobolev imbedding $|u_\lambda^{+}|^\beta\in L^{2q'}(\R^2)$. Thus,
$$
|w_{L,\lambda}^{+}|_\gamma^2 \leq C\beta^2\left(\int_{\R^2}|u_\lambda^{+}(u_{L,\lambda}^{+})^{\beta-1}|^{2q'}\right)^{1/q'}\leq C\beta^2\left(\int_{\R^2}|u_\lambda^{+}|^{2q'\beta}\right)^{1/q'}<+\infty
$$
Using Fatous' lemma in the variable $L$, we derive
$$
\left(\int_{\R^2}|u_\lambda^{+}|^{\gamma\beta}\right)^{2q'/\gamma}\leq C^{2q'}\beta^{2q'}\int_{\R^2}|u_\lambda^{+}|^{2q'\beta},
$$
from where it follows that
\begin{equation}\label{u13}
|u_\lambda^{+}|_{\beta\gamma}\leq C^{1/\beta}\beta^{1/\beta}|u_\lambda^{+}|_{\beta 2q'},
\end{equation}
Now, fixing $\gamma>2q'$, we get 
\begin{equation}\label{low2}
|u_\lambda^{+}|_\infty\leq \bar{K},~~~~ \forall \lambda\geq 1.
\end{equation}
Analogously, if we define for each $\lambda\geq 1$, $L>0$ and $\beta>1$, the functions $u_\lambda^{-}=\max\{-u_\lambda,0\}$,
$$
u_{L,\lambda}^{-}:=\left\{
\begin{array}{lll}
u_{\lambda}^{-},& \mbox{if}~~~u_\lambda\geq -L,\\[0.3cm]
L,&\mbox{if}~~~u_\lambda\leq -L,
\end{array}
\right.
$$
$$
z_{L,\lambda}^{-}:=u_\lambda^{-}(u_{L,\lambda}^{-})^{2(\beta-1)}~~~~~~ \mbox{and}~~~~~~ w_{L,\lambda,i}^{-}:=u_\lambda^{-}(u_{L,\lambda}^{-})^{\beta-1},
$$
we can prove that
\begin{equation}\label{low1}
|u_\lambda^{-}|_{\infty}\leq \tilde{K},~~~
\forall\lambda\geq 1.
\end{equation}
Therefore, from $(\ref{low2})$ and $(\ref{low1})$, 
\begin{equation}\label{es223}
|u_\lambda|_{\infty}\leq K,~~ \forall\lambda\geq 1,
\end{equation}
for some $K>0$, which proves the proposition.

\fim

\begin{proposition}\label{p5}
Let $\{ u_\lambda\}\subset\mathcal{H}_\lambda$ be a family of nodal solution of $(A)_\lambda$
with $\|u_\lambda\|^2\leq m<1$ and $u_\lambda\rightarrow 0$ in $H^1(\R^2\setminus\Omega_\Gamma)$
as $\lambda\rightarrow \infty$. Then, there exists $\lambda^{*}>0$
with the following property:
$$
|u_\lambda|_{\infty,\R^2\setminus\Omega_\Gamma'}\leq
a,~~~\forall\lambda\geq\lambda^{*}.
$$
Hence, $u_\lambda$ is a nodal solution of $(P)_\lambda$ for
$\lambda\geq\lambda^{*}$.
\end{proposition}

\noindent {\bf Proof.} In this proof, we adapt some arguments explored in \cite{AF} to get an estimate for the $L^{\infty}$-norm of the family $\{u_\lambda\}$ on a neighborhood of $\partial\Omega_\Gamma'$. In doing so, we will conclude easily the proof of the proposition.

Let $x_1,...,x_l\in\partial\Omega_\Gamma'$, $R>0$ and $0<r<R/2$ such that
$$
\partial\Omega_\Gamma'\subset\mathcal{N}(\partial\Omega_\Gamma'):=\bigcup_{i=1}^lB_{R+r}(x_i)
$$
and
$$
B_{R+r}(x_i)\subset\R^2\setminus\Omega_\Gamma,~~~\forall
i\in\{1,...,l\}.
$$

Consider $\eta_i\in C^\infty(\R^2)$, $0\leq\eta_i\leq 1$ with
$$
\eta_i(x)=\left\{\begin{array}{lll}
1,&~~\mbox{if}~~~~|x-x_i|\leq R\\[0.3cm]
0,&~~\mbox{if}~~~~|x-x_i|\geq R+r
\end{array}
\right.
$$
and $|\nabla \eta_i|\leq 2/r$, for each $i\in\{1,...,l\}$.

Now, for each $\lambda\geq 1$, $L>0$ and $\beta>1$, let us define
$$
u_{L,\lambda}^{+}:=\left\{
\begin{array}{lll}
u_{\lambda}^{+},& \mbox{if}~~~u_\lambda\leq L,\\[0.3cm]
L,&\mbox{if}~~~u_\lambda\geq L,
\end{array}
\right.
$$
$$
z_{L,\lambda,i}^{+}:=\eta_i^2u_\lambda^{+}(u_{L,\lambda}^{+})^{2(\beta-1)}~~~~~~ \mbox{and}~~~~~~ w_{L,\lambda,i}^{+}:=\eta_i u_\lambda^{+}(u_{L,\lambda}^{+})^{\beta-1}.
$$

Repeating the same ideas employed in the proof of Proposition \ref{p5}, it follows that

\begin{equation}\label{zu14}
|w_{L,\lambda,i}^{+}|_\gamma^2\leq C\beta^2\left[ \int_{\R^2}|u_\lambda^{+}|^2(u_{L,\lambda}^{+})^{2(\beta-1)}|\nabla\eta_i|^2+\int_{\R^2}b_\tau(u_\lambda^{+})\eta_i^2|u_\lambda^{+}|^2(u_{L,\lambda}^{+})^{2(\beta-1)}\right].
\end{equation}

Using Proposition $\ref{p5}$ and the definition of $b_\tau$, we obtain $|b_\tau(u_\lambda^{+})|_\infty\leq C$, for all $\lambda\geq 1$ and some constant $C>0$. Then,  from definition of $\eta_i$ and $(\ref{zu14})$,
$$
\left(\int_{B_R(x_i)}|u_\lambda^{+}|^{\gamma}(u_{L,\lambda})^{\gamma(\beta-1)}\right)^{2/\gamma}\leq C\beta^2\int_{B_{R+r(x_i)}}|u_\lambda^{+}|^{2\beta}.
$$

Using Fatous' lemma in the variable $L$, we obtain
\begin{equation}\label{zu15}
\left(\int_{B_R(x_i)}|u_\lambda^{+}|^{\gamma\beta}\right)^{2/\gamma}\leq C\beta^2\int_{B_{R+r(x_i)}}|u_\lambda^{+}|^{2\beta}.
\end{equation}

Now, if $\beta=\dfrac{\gamma(t-1)}{2t}$ with $t=\dfrac{\gamma^2}{2(\gamma-2)}$, then $\beta>1$, $\dfrac{2t}{t-1}<\gamma$ and $u_\lambda^{+}\in L^{\beta 2t/(t-1)}(B_{R+r}(x_i))$. It follows from $(\ref{zu15})$ and H\"older's inequality with exponents $t/(t-1)$ and $t$ that
$$
\left(\int_{B_R(x_i)}|u_\lambda^{+}|^{\gamma\beta}\right)^{2/\gamma}\leq C\beta^2\left[\int_{B_{R+r}(x_i)}|u_\lambda^{+}|^{2\beta t/(t-1)}\right]^{(t-1)/t}\left[\int_{B_{R+r}(x_i)}1\right]^{1/t},
$$
that is,
\begin{equation}\label{zu16}
|u_\lambda^{+}|_{L^{\gamma\beta}(B_R(x_i))}\leq C^{1/\beta}\beta^{1/\beta}|u_\lambda^{+}|_{L^{2\beta t/(t-1)}(B_{R+r}(x_i))}.
\end{equation}

If we consider $\chi=\dfrac{\gamma(t-1)}{2t}$ and $s=\dfrac{2t}{t-1}$, the inequality in $(\ref{zu16})$ gives 

\begin{equation}\label{zu17}
|u_\lambda^{+}|_{L^{\chi^{n+1}s}(B_R(x_i))}\leq C^{\sum_{i=1}^{n}\chi^{-i}}\chi^{\sum_{i=1}^{n}i\chi^{-i}}|u_\lambda^{+}|_{L^{\gamma}(B_{R+r}(x_i))},
\end{equation}
implying that 
$$
|u_\lambda^{+}|_{L^\infty(B_R(x_i))}\leq C |u_\lambda^{+}|_{L^{\gamma}(B_{R+r}(x_i))}.
$$ 
Using the convergence of $(u_\lambda^{+})$ to $0$ in $H^1(\R^2\setminus\Omega_\Gamma)$ as $\lambda\to\infty$, for a fixed $\epsilon>0$, there exists $\lambda_{\epsilon,i}\geq 1$ such that
$$
|u_\lambda^{+}|_{L^\infty(B_R(x_i))}\leq \epsilon,~~~~ \forall\lambda\geq\lambda_{\epsilon,i}.
$$
In particular, fixing $\epsilon=a$ and $\lambda_{*}=\displaystyle\max_{1\leq i\leq l}\{\lambda_{a,i}\}$, we conclude that
\begin{equation}\label{es2}
|u_\lambda^{+}|_{\infty,\mathcal{N}(\partial\Omega_\Gamma')}\leq a
~~~~\mbox{for all}~~~~ \lambda\geq\lambda^{*}.
\end{equation}
Analogously, if we define for each $\lambda\geq 1$, $L>0$ and $\beta>1$, the functions $u_\lambda^{-}=\max\{-u_\lambda,0\}$,
$$
u_{L,\lambda}^{-}:=\left\{
\begin{array}{lll}
u_{\lambda}^{-},& \mbox{if}~~~u_\lambda\geq -L,\\[0.3cm]
L,&\mbox{if}~~~u_\lambda\leq -L,
\end{array}
\right.
$$
$$
z_{L,\lambda,i}^{-}:=\eta_i^2u_\lambda^{-}(u_{L,\lambda}^{-})^{2(\beta-1)}~~~~~~ \mbox{and}~~~~~~ w_{L,\lambda,i}^{-}:=\eta_i u_\lambda^{-}(u_{L,\lambda}^{-})^{\beta-1},
$$
we can prove that
\begin{equation}\label{es222}
|u_\lambda^{-}|_{\infty,\mathcal{N}(\partial\Omega_\Gamma')}\leq a
~~~~\mbox{for all}~~~~ \lambda\geq\lambda^{*}.
\end{equation}
Thus, from $(\ref{es2})$ and $(\ref{es222})$, 
\begin{equation}\label{es223}
|u_\lambda|_{\infty,\mathcal{N}(\partial\Omega_\Gamma')}\leq a
~~~~\mbox{for all}~~~~ \lambda\geq\lambda^{*}.
\end{equation}

Now, for $\lambda\geq\lambda^{*}$, we define
$v_\lambda:\R^2\setminus\Omega_\Gamma'\to\R$ by
$$
v_\lambda(x)=\left(u_\lambda(x)-a\right)^{+}.
$$
From $(\ref{es2})$, we have $v_\lambda\in H_0^1(\R^2\setminus\Omega_\Gamma')$. Using $v_\lambda$ as a test function, it is possible to prove that $v_\lambda\equiv 0$ in $\R^2\setminus\Omega_\Gamma'$, that is,  $u_\lambda(x)\leq a$ a.e. in $\R^2\setminus\Omega_\Gamma'$. Considering the function $(u_\lambda-a)^{-}$, the same argument works well to show that $u_\lambda(x)\geq -a$ a.e. in $\R^2\setminus\Omega_\Gamma'$. Thus, $|u_\lambda(x)|\leq a$ a.e. in $\R^2\setminus\Omega_\Gamma'$. Therefore, from Remark \ref{rea}, the proof is finished.  

\fim

\section{A special family of nodal solution to $(A)_\lambda$}\label{scv}

In this section, as in \cite{A1}, we modify all the sets that appear in the minimax arguments explored in \cite{DT} to get nodal solutions. These modifications are necessary, because we are working with exponential critical growth, and in this case, the estimate involving the norm of sequences must be very carefully obtained to use the Tundiger-Moser inequalities mentioned in the introduction of a correct way. After that, using deformation lemma, we show the existence of a special family of nodal solutions to $(A)_\lambda$ for $\lambda$ large enough. These nodal solutions are exactly the nodal solutions given in Theorem \ref{bd}.

In what follows, let us fix $\epsilon>0$ and $\zeta=\zeta(\epsilon)>0$ such that
\begin{equation}\label{9}
I_j((1-\epsilon)w_j^{\pm}),~~
I_j((1+\epsilon)w_j^{\pm})<I_j(w_j^{\pm})-\zeta,~~\mbox{for
all}~~j\in\Gamma.
\end{equation}
Moreover, without loss of generality, we can assume
$\Gamma=\{1,...,l\}~(l\leq k)$. In the sequel, we denote by
$Q=(1-\epsilon,1+\epsilon)^{2l}$ and define $\gamma_0:\overline{Q}\to
\mathcal{H}_\lambda$ by
\begin{equation}\label{10}
\gamma_0(\overrightarrow{s},\overrightarrow{t})(x)=\sum_{j=1}^ls_jw_j^{+}+\sum_{j=1}^lt_jw_j^{-},
\end{equation}
where $(\overrightarrow{s},\overrightarrow{t})=(s_1,...,s_l,t_1,...,t_l)$, and the number
$$
S_{\lambda,\Gamma}=\inf_{\gamma\in\sum_\lambda}~\max_{(\overrightarrow{s},\overrightarrow{t})\in\overline{Q}}\Phi_\lambda(\gamma(\overrightarrow{s},\overrightarrow{t})),
$$
where
$$
\Sigma_{\lambda}=\left\{\gamma\in
C(\overline{Q},\mathcal{H}_\lambda):~\gamma^{\pm}|_{\Omega_j'}\neq 0,~\forall
j\in\Gamma~\mbox{and}~(\overrightarrow{s},~\overrightarrow{t})\in
\overline{Q},\gamma=\gamma_0~\mbox{on}~\partial Q\right\}.
$$
We remark that $\gamma_0\in\Sigma_\lambda$, so
$\Sigma_\lambda\neq\emptyset$ and $S_{\lambda,\Gamma}$ is well
defined.
\begin{lemma}\label{l3}
For any $\gamma\in\Sigma_\lambda$ there exists
$(\overrightarrow{s}_{*},\overrightarrow{t}_{*})\in \overline{Q}$ such that
$$
\Phi_{\lambda,j}'\left(\gamma^{\pm}(\overrightarrow{s}_{*},\overrightarrow{t}_{*})\right)\left(\gamma^{\pm}(\overrightarrow{s}_{*},\overrightarrow{t}_{*})\right)=0
$$
for all $j\in\{1,...,l\}$. As an immediate consequence,
$$
\Phi_{\lambda,j}\left(\gamma^{\pm}(\overrightarrow{s}_{*},\overrightarrow{t}_{*})\right) \geq d_{\lambda,j}.
$$
\end{lemma}
\noindent{\bf Proof.} The proof follows as in \cite[Lemma 4.1]{A1}, because the growth of $f$ is not relevant in this lemma. \fim

From now on, we denote by $D_{\Gamma}$ the number
$D_{\Gamma}=\displaystyle\sum_{j=1}^{l}d_j.$

\begin{proposition}\label{p2}
The numbers $D_\Gamma$ and $S_{\lambda,\Gamma}$ verify the following
relations
\begin{enumerate}
\item[$(a)$]$\displaystyle\sum_{j=1}^ld_{\lambda,j}\leq
S_{\lambda,\Gamma}\leq D_{\Gamma}$ for all $\lambda\geq 1$ and;
\item[$(b)$] $S_{\lambda,\Gamma}\rightarrow D_{\Gamma}$ as
$\lambda\rightarrow\infty$.
\end{enumerate}
\end{proposition}
\noindent{\bf Proof.}

$(a)$ Since $\gamma_0$ defined in $(\ref{10})$ belongs to
$\Sigma_\lambda$, 
$$
\begin{array}{rl}
S_{\lambda,\Gamma}&\displaystyle\leq\max_{(\overrightarrow{s},\overrightarrow{t})\in
\overline{Q}}\Phi_\lambda(\gamma_0(\overrightarrow{s},\overrightarrow{t}))\\[0.4cm]
&\leq\displaystyle\max_{\overrightarrow{s}\in[1-\epsilon,1+\epsilon]^l}\sum_{j=1}^{l}I_j(s_jw_j^{+})+\max_{\overrightarrow{t}\in[1-\epsilon,1+\epsilon]^l}\sum_{j=1}^{l}I_j(t_jw_j^{-}).
\end{array}
$$
From the definition of $w_j$, it is well known that
\begin{equation}\label{equ2}
\max_{z\in[1-\epsilon,1+\epsilon]}I_j(zw_j^{\pm})=I_j(w_j^{\pm}),~~~~\mbox{for
each}~~~j\in\Gamma,
\end{equation}
and thus
$$
S_{\lambda,\Gamma}\leq\sum_{j=1}^{l}d_j=D_\Gamma.
$$
Now, for $\gamma\in\Sigma_\lambda$, let
$(\overrightarrow{s}_{*},\overrightarrow{t}_{*})\in Q$ given by
Lemma \ref{l3}. Recalling that
$\Phi_{\lambda,\R^2\setminus\Omega_\Gamma'}(u)\geq 0$ for all $u\in
H^1(\R^2\setminus\Omega_\Gamma')$, we have
$$
\Phi_\lambda(\gamma(\overrightarrow{s}_{*},\overrightarrow{t}_{*}))\geq\sum_{j=1}^l\Phi_{\lambda,j}(\gamma(\overrightarrow{s}_{*},\overrightarrow{t}_{*})) \geq \sum_{j=1}^ld_{\lambda,j},
$$
and so,
$$
\max_{(\overrightarrow{s},\overrightarrow{t})\in
\overline{Q}}\Phi_\lambda(\gamma(\overrightarrow{s},\overrightarrow{t}))\geq\sum_{j=1}^ld_{\lambda,j}.
$$
Thereby, from definition of $S_{\lambda,\Gamma}$, 
$$
S_{\lambda,\Gamma}\geq\sum_{j=1}^ld_{\lambda,j},
$$
finishing the proof of $(a)$.

$(b)$ We begin proving that $d_{\lambda,j}\to d_j$ as $\lambda\to
\infty$. In fact, let $(\lambda_n)$ be an arbitrary sequence with
$\lambda_n \to +\infty$. Now, let
$w_{\lambda_n,j}\in H^1(\Omega_j')$ be a least energy nodal solution
to Problem $(\ref{3})$ given in Section \ref{snd}, with
$\lambda=\lambda_n$, i.e.
\begin{equation}\label{a1}
\Phi_{\lambda_n,j}(w_{\lambda_n,j})=d_{\lambda_n,j}~~~\mbox{and}~~~\Phi_{\lambda_n,j}'(w_{\lambda_n,j})=0
\end{equation}
The same arguments used in proof of Proposition \ref{p1} work to
prove that, for each $j\in\Gamma$ and for a subsequence
$(w_{\lambda_{n_k},j})$, there exists $w_{0,j}$ such that
$$
w_{\lambda_{n_k},j}\to w_{0,j}~~~\mbox{in}~~~
H^1(\Omega_j')~~~~\mbox{as}~~~~n_k\to\infty.
$$
Furthermore, $w_{0,j}\in H_0^1(\Omega_j)$ is a nodal
solution of Problem $(\ref{2})$. Therefore, 
$$
\lim_{k\to\infty}\Phi_{\lambda_{n_k},j}(w_{\lambda_{n_k},j})=I_j(w_{0,j}) \geq d_j.
$$
Once $d_{\lambda,j}\leq d_j$, we conclude that $d_{\lambda,j}\to d_j$ as $\lambda\to\infty$, from where
it follows that
$$
\sum_{j=1}^{l}d_{\lambda,j}\to
D_\Gamma,~~~\mbox{as}~~~\lambda\to\infty.
$$
The last limit together with $(a)$ implies that $(b)$ holds.
 \fim

Hereafter, $E_{\lambda,j}^{+}$ and $E_{\lambda,j}^{-}$ denote the
cone of nonnegative and non-positive functions belong to
$\mathcal{H}_\lambda(\Omega_j')$, respectively, that is,
$$
E_{\lambda,j}^{+}=\left\{\ u\in\mathcal{H}_\lambda(\Omega_j'):~
u(x)\geq 0~\mbox{a.e. in}~\in \Omega_j'\right\}
$$
and
$$
E_{\lambda,j}^{-}=\left\{\ u\in\mathcal{H}_\lambda(\Omega_j'):~
u(x)\leq 0~\mbox{a.e. in}~\in \Omega_j'\right\}.
$$
From the definition of $\gamma_0$, there exists a positive constant
$\tau$ such that
$$
dist_{\lambda,j}\left(\gamma_0(\overrightarrow{s},\overrightarrow{t})|_{\Omega_j'},E_{\lambda,j}^{\pm}\right)>\tau~~\mbox{for
all}~~(\overrightarrow{s},\overrightarrow{t})\in Q,~j\in
\Gamma~~\mbox{and}~~\lambda>0,
$$
where $dist_{\lambda,j}\left(K,F\right)$ denotes the distance
between sets of $\mathcal{H}_\lambda(\Omega_j')$. Taking the number
$\tau$ obtained in the last inequality, we define
$$
\Theta=\left\{
u\in\mathcal{H}_\lambda:~dist_{\lambda,j}\left(u|_{\Omega_j'},
E_{\lambda,j}^{\pm}\right)\geq\tau~~\forall j\in \Gamma\right\}.
$$
Moreover, for any $c,\mu>0$ and $0<\delta<\tau/2$, we set the
sets
$$
\Phi_\lambda^c=\left\{u\in\mathcal{H}_\lambda:~
\Phi_\lambda(u)\leq c\right\}~~\mbox{and}~~
B_{\lambda,\mu}=\left\{u\in\Theta_{2\delta}:~~
\left|\Phi_\lambda(u)-S_{\lambda,j}\right|\leq\mu\right\},
$$
where $\Theta_r$, for $r>0$, denotes the set
$\Theta_r=\left\{u\in\mathcal{H}_\lambda:~
dist_{\lambda,j}\left(u,\Theta\right)\leq r\right\}$.

Notice that for each $\mu>0$, there exists
$\Lambda^{*}=\Lambda^{*}(\mu)>0$ such that
$$
w=\sum_{j=1}^{l}w_j\in B_{\lambda,\mu},~~\mbox{for
all}~~\lambda\geq\Lambda^{*},
$$
because $w\in\Theta,~ \Phi_\lambda(w)=D_\Gamma$ and
$S_{\lambda,\Gamma}\rightarrow D_\Gamma$ as
$\lambda\rightarrow\infty$. Therefore,
$B_{\lambda,\mu}\neq\emptyset$ for $\lambda$ sufficiently large.

Observe that, for $\epsilon>0$ small enough,
$$
\left\|\gamma_0(\vec{s},\vec{t})\right\|_\lambda^2\leq (1+\epsilon)^2\sum_{j=1}^{l}\|w_j\|_{\Omega_j}^2\leq M:=\dfrac{2\theta D_\Gamma}{\theta -2}(1+\epsilon)^2<1.
$$
The fact that $M<1$ is crucial in our argument, because we are working with exponential critical growth, see for example, Claim \ref{cl1} below . However, this type of analysis is not necessary when $N \geq 3$, see  \cite{A2}, \cite{A1} and \cite{AMFS}.

In the sequel, for $r>0$, let us consider 
$$
\overline{B}_{r}(0)=\left\{
u\in\mathcal{H}_\lambda:~\|u\|_\lambda\leq r\right\},
$$ 
and we denote by $\mu^{*}$ the ensuing real number
\begin{equation}\label{11}
\mu^{*}=\min\left\{\dfrac{M+1}{2},\dfrac{\delta}{2}\right\}.
\end{equation}

As a consequence of the above consideration, we have the following result.

\begin{proposition}\label{p3}
For each $\mu>0$ fixed, there exist $\sigma_o=\sigma_o(\mu)>0$ and
$\Lambda_{*}=\Lambda_{*}(\mu)\geq 1$ independent of $\lambda$ such
that
$$
\|\Phi_\lambda'(u)\|_\lambda^{*}\geq
\sigma_o~~\mbox{for}~~\lambda\geq\Lambda_{*}~~\mbox{and all}~~u\in
\left(B_{\lambda,2\mu}\setminus
B_{\lambda,\mu}\right)\cap\overline{B}_{(M+3)/4}(0)\cap
\Phi_\lambda^{D_\Gamma}.
$$
\end{proposition}
\noindent{\bf Proof.} Arguing by contradition, we assume that there exist $\lambda_n\to\infty$ and
\begin{equation}\label{doi1.1}
u_n\in(B_{\lambda_n,2\mu}\setminus B_{\lambda_n,\mu}) \cap \overline{B}_{(M+3)/4}(0)\cap\Phi^{D_\Gamma}
\end{equation}
such that $\|\Phi_{\lambda_n}'(u_n)\|_{\lambda_n}^{*}\to 0$, as $\lambda_n\to\infty$. Since  
$$
S_{\lambda_n,\Gamma}-\Phi_{\lambda_n}(u_n)\leq|\Phi_{\lambda_n}(u_n)-S_{\lambda_n,\Gamma}|\leq 2\mu,
$$
and  $S_{\lambda_n,\Gamma}= D_\Gamma+o_n(1)$, we derive
$$
D_\Gamma-2\mu+o_n(1)\leq\Phi_{\lambda_n}(u_n)\leq D_\Gamma.
$$
Thus, $(\Phi_{\lambda_n}(u_n))$ is a bounded sequence and we may suppose
$$
\Phi_{\lambda_n}(u_n)\to d\in[D_\Gamma-2\mu, D_\Gamma],
$$
after extracting a subsequence. Applying Proposition \ref{p1}, we can extract a subsequence $u_n\to u$ in $H^1(\R^2)$, where $u$ is a solution of $(P)_j$ with
$$
\|u_n-u\|_{\lambda_n}\to 0,~~\lambda_n\int_{\R^2}V(x)|u_n|^2\to 0~~\mbox{and}~~\|u_n\|_{\lambda_n,\R^2\setminus\Omega_\Gamma}\to 0.
$$
Once $u_n\in \Theta_{2\delta}$, we have that $\|u_n^{\pm}\|_{\lambda_n,\Omega_j'}\geq\tau-2\delta>0$, for all $\lambda_n$, leading to
$\|u^{\pm}\|_{\Omega_j}\neq 0$, for all $j\in\Gamma$. Then, $u$ is a nodal solution of $(P)_j$, for all $j\in \Gamma$, and
$$
\sum_{j=1}^{l}d_j\leq \sum_{j=1}^lI_j(u|_{\Omega_j})\leq D_\Gamma.
$$
This fact gives $I_j(u|_{\Omega_j})=d_j$, for all $j\in \Gamma$, and hence $\Phi_{\lambda_n}(u_n)\to D_\Gamma$. On the other hand, since $S_{\lambda_n,\Gamma}\to D_\Gamma$, we have
$$
|\Phi_{\lambda_n}(u_n)-S_{\lambda_n,\Gamma}|\to 0,~~\mbox{as}~~\lambda_n\to\infty.
$$ 
Therefore, $u_n\in B_{\lambda_n,\mu}$, for $n$ large enough, which contradicts (\ref{doi1.1}).
\fim
\begin{proposition}\label{doi2}
For each $\mu\in(0,\mu^{*})$, there exists $\Lambda^{*}=\Lambda^{*}(\mu)>0$ such that for all $\lambda\geq\Lambda^{*}$ the functional $\Phi_\lambda$ has a critical point in $B_{\lambda,\mu}\cap\overline{B}_{(M+3)/4}\cap\Phi_{\lambda}^{D_\Gamma}$.

\end{proposition}
\noindent{\bf Proof.} Arguing again by contradiction, we suppose that there exist $\mu\in(0,\mu^{*})$ and a sequence $\lambda_n\to+\infty$, such that $\Phi_{\lambda_n}$ has no critical points in $B_{\lambda_n,\mu}\cap\overline{B}_{(M+3)/4}\cap\Phi_{\lambda_n}^{D_\Gamma}$. From Proposition \ref{psc}, the $(PS)_c$ condition holds for $\Phi_{\lambda_n}$, for $c\in(0,D]$. Thus, there exists a constant $d_{\lambda_n}>0$ such that
$$
\|\Phi_{\lambda_n}'(u_n)\|_{\lambda_n}^{*}\geq d_{\lambda_n}~~\mbox{for all}~~u\in B_{\lambda_n,\mu}\cap\overline{B}_{(M+3)/4}(0)\cap\Phi_{\lambda_n}^{D_\Gamma}.
$$
Moreover, from Proposition \ref{p3}, we also have
$$
\|\Phi_{\lambda_n}'(u)\|_{\lambda_n}^{*}\geq
\sigma_o~~\mbox{for all}~~u\in
\left(B_{\lambda_n,2\mu}\setminus
B_{\lambda_n,\mu}\right)\cap\overline{B}_{(M+3)/4}(0)\cap
\Phi_{\lambda_n}^{D_\Gamma}
$$
and for all $\lambda_n\geq\Lambda_{*}$, where $\sigma_o>0$ is independent of $\lambda_n$, for $n$ large enough.

In what follows, $\Psi_n:\mathcal{H}_{\lambda_n}\to\R~\mbox{and}~H_n:\Phi_{\lambda_n}^{D_{\Gamma}}\to\mathcal{H}_{\lambda_n}$
are continuous functions verifying
$$
\begin{array}{rl}
\Psi_n(u)=1,&\mbox{for}~~u\in B_{\lambda_n,3\mu/2}\cap\Theta_\delta\cap\overline{B}_{(M+1)/2}(0),\\[0.2cm]
\Psi_n(u)=0,&\mbox{for}~~u\notin B_{\lambda_n,3\mu/2}\cap\Theta_\delta\cap\overline{B}_{(M+3)/4}(0),\\[0.2cm]
0\leq\Psi_n(u)\leq 1,&\mbox{for}~~ u\in\mathcal{H}_{\lambda_n},
\end{array}
$$
and
$$
H_n(u)=\left\{
\begin{array}{rcl}
-\Psi_n(u)\|Y_n(u)\|^{-1}Y_n(u),&\mbox{for}~~ u\in B_{\lambda_n,2\mu}\cap \overline{B}_{(M+3)/4}(0),\\[0.2cm]
0,&\mbox{for}~~u\notin B_{\lambda_n,2\mu}\cap \overline{B}_{(M+3)/4}(0),
\end{array}
\right.
$$
where $Y_n$ is a pseudo-gradient vector field for $\Phi_{\lambda_n}$ on 
$$
\mathcal{M}_n=\{u\in\mathcal{H}_{\lambda_n}:~\Phi_{\lambda_n}'(u)\neq 0\}.
$$
From the definition of $H_n$, 
$$
\|H_n(u)\|\leq 1,~~\mbox{for all}~~u\in\Phi_{\lambda_n}^{D_\Gamma}.
$$
Hence, there exists a deformation flow $\eta_n:[0,\infty)\times \Phi_{\lambda_n}^{D_\Gamma}\to\Phi_{\lambda_n}^{D_\Gamma}$ given by
$$
\dfrac{d\eta_n}{dt}=H_n(\eta_n),~~~~\eta_n(0,u)=u\in \Phi_{\lambda_n}^{D_\Gamma}.
$$
This flow satisfies the following basic properties:
\begin{equation}\label{v1}
\dfrac{d}{dt}\Phi_{\lambda_n}(\eta_n(t,u))\leq -\Psi_n(\eta_n(t,u))\|\Phi_{\lambda_n}'(\eta_n(t,u))\|\leq 0,
\end{equation}

\begin{equation}\label{v2}
\left\|\dfrac{d\eta_n}{dt}\right\|_{\lambda_n}=\|H_n(\eta_n)\|_{\lambda_n}\leq 1
\end{equation}
and
\begin{equation}\label{v3}
\eta_n(t,u)=u,~~~\forall t\geq 0,~~ u\notin B_{\lambda_n,2\mu}\cap \overline{B}_{(M+3)/4}(0).
\end{equation}

Next, let us show that the functions $\gamma_n:\overline{Q}\to\mathcal{H}_{\lambda_n}$ belongs to $\Sigma_{\lambda_n}$, for $n$ large enough. We begin observing that $\gamma_n$ is a continuous function in $\overline{Q}$. Since $\mu\in(0,\mu^{*})$, from $(\ref{9})$, $(\ref{equ2})$ and $(\ref{11})$, 
$$
|\Phi_{\lambda_n}(\gamma_0(\vec{s},\vec{t}))-D_\Gamma|>\zeta\geq\delta\geq 2\mu^{*},~~~\forall (\vec{s},\vec{t})\in\partial Q,~n\in\mathbb{N}.
$$
Thus, using again the fact that $S_{\lambda,\Gamma}\to D_\Gamma$ quando $\lambda\to\infty$, there exists $n_o>0$ such that
$$
|\Phi_{\lambda_n}(\gamma_0(\vec{s},\vec{t}))-S_{\lambda_n,\Gamma}|>2\mu,~\forall (\vec{s},\vec{t})\in\partial Q,~n\geq n_o,
$$
which implies that $\gamma_0(\vec{s},\vec{t})\notin B_{\lambda_n,2\mu}$, for all $(\vec{s},\vec{t})\in\partial Q$ and $n\geq n_o$. So,
$$
\eta_n(t,\gamma_0(\vec{s},\vec{t}))=\gamma_0(\vec{s},\vec{t})~ \mbox{for all}~ (\vec{s},\vec{t})\in\partial Q.
$$
Now, we only have to prove that
$$
\gamma_n(\vec{s},\vec{t})^{\pm}\in H^1(\Omega_j')\setminus\{0\},
$$
for all $j\in\Gamma$ and $(\vec{s},\vec{t})\in Q$.

Once that $\gamma_n(\vec{s},\vec{t})=\eta_n(T_n,\gamma_o(\vec{s},\vec{t}))\in\Theta_{2\delta}$ for all $n$, we have
$$
dist_{\lambda_n,j}(\gamma_n(\vec{s},\vec{t}), E_{\lambda_n,j}^{\pm})\geq \tau-2\delta>0.
$$
Then, $\gamma_n^{\pm}|_{\Omega_j}\neq 0$ for all $j\in\Gamma$, implying that $\gamma_n\in\Sigma_{\lambda_n}$ for $n$ large enough.

Note that $supt~\gamma_0(\vec{s},\vec{t})\subset\overline{\Omega}_{\Gamma}$ for all $(\vec{s},\vec{t})\in\overline{Q}$ and that $\Phi_{\lambda}(\gamma_0(\vec{s},\vec{t}))$ does not depend on $\lambda\geq 1$. Furthermore, 
$$
\Phi_{\lambda}(\gamma_0(\vec{s},\vec{t}))\leq D_\Gamma,~~\mbox{for all}~~(\vec{s},\vec{t})\in \overline{Q}
$$
and
$$
\Phi_{\lambda}(\gamma_0(\vec{s},\vec{t}))= D_\Gamma~~\mbox{if, and only if,}~~s_j=t_j=1,~~\forall j\in\{1,...,l\}.
$$
Therefore, the number
$$
m_0^{n}:=\sup\left\{\Phi_{\lambda_n}(u):~u\in\gamma_0(Q)\setminus(B_{\lambda_n,\mu}\cap\overline{B}_{\frac{M+1}{2}}(0))\right\},
$$
is independent of $\lambda_n$ and verifies
$$
\limsup_{n\to\infty}m_0^n<D_\Gamma.
$$

The next claim is crucial, because we are working with exponential growth, and some arguments used in \cite{A2}, \cite{A1} and \cite{AMFS}  cannot be used directly, so a careful analysis is necessary.

\begin{claim}\label{cl1}
There exists a constant $K>0$ such that
$$
|\Phi_{\lambda_n}(u)-\Phi_{\lambda_n}(v)|\leq K\|u-v\|_{\lambda_n}
$$
for all $u,v\in \overline{B}_{(M+3)/4}(0)$.
\end{claim}
In fact, let $u,v\in\overline{B}_{(M+3)/4}(0)$, there is $K>0$ such that
$$
|\langle\Phi_{\lambda_n}'(tu+(1-t)v),w\rangle|\leq K,~~~\forall w\in \mathcal{H}_{\lambda_n},~\|w\|_{\lambda_n}\leq 1.
$$
Since,
$$
|\langle\Phi_{\lambda_n}'(tu+(1-t)v),w\rangle|\leq \dfrac{M+3}{2}+\int_{\R^2}|g(tu+(1-t)v)w|,
$$
we only need to prove the boundedness of the above integral. Using the growth of $g$ given in $(\ref{g})$,
\begin{equation}\label{equa1}
\int_{\R^2}|g(tu+(1-t)v)w|\leq\dfrac{M+3}{2}+C\int_{\R^2}|w|b_\tau(tu+(1-t)v).
\end{equation}
By the H\"older's inequality, 
\begin{equation}\label{equa2}
\int_{\R^2}|w|b_\tau(tu+(1-t)v)\leq |w|_{q'}|b_\tau(tu+(1-t)v)|_q,
\end{equation}
where $1/q+1/q'=1$. Since $M<1$,
$$
\|tu+(1-t)v\|_{\lambda_n}\leq t\|u\|_{\lambda_n}+(1-t)\|v\|_{\lambda_n}\leq\dfrac{M+3}{4}<1.
$$
Then, we can take $q>1$, $q$ near $1$, such that $
q\tau(M+3)/4<1$. Thus, from Corollary \ref{ctec}
\begin{equation}\label{equa3}
|b_\tau(tu+(1-t)v)|_q\leq C,~~\forall t\in[0,1],~u,v\in B_{(M+3)/4}(0).
\end{equation}
Therefore, from $(\ref{equa1})$, $(\ref{equa2})$ and $(\ref{equa3})$, 
$$
\int_{\R^2}|f(tu+(1-t)v)w|\leq C,~~~~\forall u,v \in \overline{B}_{(M+3)/4}(0),~ t\in[0,1],~\|w\|_\lambda\leq 1, 
$$
showing that Claim \ref{cl1} holds.

\vspace{0.5 cm}

As a consequence of the above considerations, we are able to repeat the arguments found in \cite{A1} to prove the ensuing claim

\begin{claim}\label{cl2}
There exists $T_n=T(\lambda_n)>0$ and $\epsilon^{*}>0$ independent of $n$ such that
$$
\limsup_{n\to\infty}\left\{\max_{(\vec{s},\vec{t})\in Q}\Phi_{\lambda_n}(\eta_n(T_n,\gamma_0(\vec{s},\vec{t})))\right\}<D_\Gamma-\epsilon^{*}.
$$
\end{claim}
The above claim gives 
$$
\limsup_{n\to\infty}S_{\lambda_n,\Gamma}\leq D_\Gamma-\epsilon^{*},
$$
which contradicts the Proposition \ref{p2}, and the proposition follows.

\fim

From the last proposition, we have the following result.

\begin{proposition}\label{c1}
For each $\mu\in (0,\mu^{*})$ fixed, there exist
$\Lambda^{*}=\Lambda^{*}(\mu)>1$ such that $(A)_\lambda$ has a nodal
solution $u_\lambda\in B_{\lambda,\mu}$ for all $\lambda\geq
\Lambda^{*}$.
\end{proposition}

\section{Proof of Theorem \ref{bd}}

From Proposition \ref{c1}, for each $\mu\in(0,\mu^{*})$ fixed, there
exists $\Lambda^{*}=\Lambda^{*}(\mu)>1$ such that the auxiliary
problem $(A)_\lambda$ has a nodal solution $u_\lambda\in
B_{\lambda,\mu}$ for $\lambda\geq\Lambda^{*}$ with
\begin{equation}\label{eq6}
dist_{\lambda,j}(u_\lambda,E^{\pm}_{\lambda,j})\geq
\tau-2\delta>0~~~~\forall j\in\Gamma.
\end{equation}
Repeating the same arguments used in the proof of Proposition
\ref{p1}, we get
$$
u_\lambda\to 0~~~\mbox{in}~~~
H^1(\R^2\setminus\Omega_\Gamma)~~~\mbox{as}~~~\lambda\to\infty.
$$
This together with Proposition \ref{p5} implies that $u_\lambda$ is
a nodal solution to $(P)_\lambda$, for $\lambda$ large enough.

Fixing $\lambda_n\to+\infty$, the sequence
$(u_{\lambda_n})$ verifies
$$
\Phi'_{\lambda_n}(u_{\lambda_n})=0~~~~\mbox{and}~~~~\Phi_{\lambda_n}(u_{\lambda_n})=S_{\lambda_n,\Gamma}+o_n(1),
$$
and so,
$$
\Phi'_{\lambda_n}(u_{\lambda_n})=0~~~~\mbox{and}~~~~\Phi_{\lambda_n}(u_{\lambda_n})=D_{\Gamma}+o_n(1).
$$
Therefore, $(u_{\lambda_n})$ is a $(PS)_{\infty,D_\Gamma}$, with
$D_\Gamma\in(0,D]$. By Proposition \ref{p1}, there exists $u\in
H_0^1(\Omega_\Gamma)$ such that, for some subsequence still denoted
by $(u_{\lambda_n})$,
$$
u_{\lambda_n}\to u~~\mbox{in}~~
H^1(\R^2),~~\lambda_n\int_{\R^2}V(x)|u_{\lambda_n}|^2\to
0~~\mbox{and}~~\|u_{\lambda_n}\|^2_{\lambda_n,\R^2\setminus\Omega_\Gamma}\to
0.
$$
From the above limits, we see that 
\begin{equation}\label{eq7}
I_j(u)=0~~~\mbox{for all}~~~ j\in
\Gamma~~~~\mbox{and}~~~~\sum_{j=1}^{l}I_j(u)=D_\Gamma.
\end{equation}
Once $(u_{\lambda_n})$ verifies $(\ref{eq6})$, we derive that 
\begin{equation}\label{ela1}
\tau-2\delta\leq \|u_{\lambda_n}^{\pm}\|_{\lambda_n,\Omega_j'},~~~\forall j\in\Gamma,~\forall n\in\mathbb{N}.
\end{equation}

Using these information, we are ready to prove the following claim 

\begin{claim}\label{clai}
There exists $\kappa_o>0$ such that
\begin{equation}\label{doi11}
\int_{\Omega_j}|u^{\pm}_{\lambda_n}|^{q'}\geq \kappa_o~~~\forall\lambda_n\geq\Lambda^{*},~\forall j\in\Gamma,
\end{equation}
for some $q'>1$.
\end{claim}
In fact, let us fix $j\in\Gamma$ and consider $\eta_i\in C^{\infty}(\R^2,\R)$ satisfying
$$
\eta_j\equiv 1~~\mbox{in}~~\Omega_j'~~\mbox{e}~~\eta_j\equiv 0~~\mbox{in}~~\R^2\setminus(\Omega_j')_\delta~~\mbox{and}~~((\Omega_j')_\delta\setminus\Omega_j')\subset\R^2\setminus\Omega_\Gamma.
$$
Taking $v_j=\eta_ju_{\lambda_n}^{+}$ as a test function, we obtain
$$
\int_{\R^2}\nabla u_{\lambda_n}\nabla (\eta_ju_{\lambda_n}^{+})+(\lambda_nV(x)+1)u_{\lambda_n}\eta_ju_{\lambda_n}^{+}=\int_{\R^2}f(u_{\lambda_n})\eta_ju_{\lambda_n}^{+},
$$
from where it follows that
\begin{equation}\label{doi0}
\int_{\R^2}|\nabla u_{\lambda_n}^{+}|^2\eta_j+(\lambda_nV(x)+1)|u_{\lambda_n}^{+}|^2\eta_j=\int_{\R^2}f(u_{\lambda_n})\eta_ju_{\lambda_n}^{+}+o_n(1).
\end{equation}
Combining the growth of $f$ and H\"older's inequality, we know that
$$
\int_{\R^2}f(u_{\lambda_n})\eta_ju_{\lambda_n}^{+}\leq\epsilon\int_{\R^2}|u_{\lambda_n}^{+}|^2\eta_j+C\left(\int_{\R^2}|\eta_ju_{\lambda_n}^{+}|^{q'}\right)^{1/q'}|b_\tau(u_{\lambda_n})|_q
$$
where $1/q+1/q'=1$. On the other hand, as in Proposition \ref{psc}, we can extract a subsequence verifying  
$$
\displaystyle\limsup_{n\to\infty}\|u_{\lambda_n}\|^{2}\leq m<1.
$$ 
Fixing $q>1$, with $q$ sufficiently close to $1$ and using Corollary \ref{ctec}, we derive
\begin{equation}\label{doi1}
\int_{\R^2}f(u_{\lambda_n})\eta_ju_{\lambda_n}^{+}\leq\epsilon\int_{\R^2}|u_{\lambda_n}^{+}|^2\eta_j+C\left(\int_{\R^2}|\eta_ju_{\lambda_n}^{+}|^{q'}\right)^{1/q'}.
\end{equation}
From (\ref{doi0}) and (\ref{doi1}),
$$
(1-\epsilon)\int_{\R^2}|\nabla u_{\lambda_n}^{+}|^2\eta_j+(\lambda_nV(x)+1)|u_{\lambda_n}^{+}|^2\eta_j\leq C\left(\int_{\R^2}|\eta_ju_{\lambda_n}^{+}|^{q'}\right)^{1/q'}+o_n(1).
$$
Thereby, fixing $\epsilon<1$ and using (\ref{ela1}), 
$$
0<(1-\epsilon)(\tau-2\delta)\leq(1-\epsilon)\|u_{\lambda_n}^{+}\|_{\lambda_n,\Omega_j'}^2\leq C\left(\int_{\Omega_j'}|u_{\lambda_n}^{+}|^{q'}\right)^{1/q'}+o_n(1),
$$
implying that there is $\kappa_o>0$,
$$
\int_{\Omega_j'}|u_{\lambda_n}^{+}|^{q'}\geq\kappa_o>0,~~~\forall \lambda_n\geq\Lambda^{*}.
$$
The same arguments work to prove that  
$$
\int_{\Omega_j'}|u_{\lambda_n}^{-}|^{q'}\geq\kappa_o>0,~~~\forall \lambda_n\geq\Lambda^{*},
$$
and the proof of Claim \ref{clai} is complete.

Passing to the limit as $n\to\infty$ in (\ref{doi11}), we derive
$$
\int_{\Omega_j'}|u^{\pm}|^{q'}=\int_{\Omega_j}|u^{\pm}|^{q'}\geq\kappa_o,~~~~\forall j\in\Gamma.
$$
Thus, $u$ changes sign in $\Omega_j$, for all $j\in\Gamma$. Consequently,
\begin{equation}\label{eq8}
I_j(u)\geq d_j,~~~~\forall j\in \Gamma.
\end{equation}
From $(\ref{eq7})$ and $(\ref{eq8})$, follows that $I_j(u)=d_j$ for all $j\in \Gamma$. This shows that, for each $j\in\Gamma$, $u|_{\Omega_j}$ is a least energy nodal solution for the problem $(\ref{2})$, finishing the proof .\fim

\end{document}